\documentclass[review]{elsarticle}
\usepackage{bm,amssymb,amsmath,mathrsfs}
\usepackage{amsthm}
\usepackage{lineno}
\usepackage[usenames]{color}

\usepackage{dcolumn}% Align table columns on decimal point

\newcommand{\bq}{\begin{equation}}
\newcommand{\eq}{\end{equation}}
\newcommand{\bc}{\begin{center}}
\newcommand{\ec}{\end{center}}
\newcommand{\bit}{\begin{itemize}}
\newcommand{\eit}{\end{itemize}}
\newcommand{\ben}{\begin{enumerate}}
\newcommand{\een}{\end{enumerate}}

\theoremstyle{plain}

\newtheorem*{theorem*}{Theorem}
\begin{document}

%\journal{\ }
%\journal{Expositiones Mathematicae}
\journal{(internal report CC23-11)}

\begin{frontmatter}

\title{Factorial moments of the Poisson distribution of order $k$}

\author[cc]{S.~R.~Mane}
\ead{srmane001@gmail.com}
\address[cc]{Convergent Computing Inc., P.~O.~Box 561, Shoreham, NY 11786, USA}

\begin{abstract}
The factorial moments of the standard Poisson distribution are well known.
The present note presents an explicit combinatorial sum for the factorial moments of the Poisson distribution of order $k$.
Unlike the standard Poisson distribution (the case $k=1$), for $k>1$ the structure of the factorial moments is much more complicated.
Some properties of the factorial moments of the Poisson distribution of order $k$ are elucidated in this note.
\end{abstract}

\vskip 0.25in

\begin{keyword}
Poisson distribution of order $k$
\sep moments
\sep factorial moments
\sep Compound Poisson distribution  
\sep discrete distribution 

%\vskip 0.25in
% MSC2020 codes here, in the form: \MSC code \sep code
\MSC[2020]{
%primary 
60E05  % Probability distributions: general theory
\sep 39B05 % General theory of functional equations and inequalities
\sep 11B37  % recurrences
\sep 05-08  % Computational methods for problems pertaining to combinatorics
}

%\vskip 0.25in
%\PACS{
%02.20.Qs % general properties, structure, and representation of Lie groups
%\sep 02.20.Hj % classical groups
%\sep 02.30.Gp % special functions
%\sep 29.27.Hj % polarized beams
%\sep 29.20.D- % cyclic accelerators and storage rings
%\sep 29.20.db % storage rings and colliders
%\sep 29.27.-a   %Beams in particle accelerators  
%\sep 41.85.-p % beam optics
%\sep 02.30.Ik %Integrable systems 
%\sep 02.60.Lj % Ordinary and partial differential equations; %boundary value problems
%\sep 13.40.Em % Electric and magnetic moments
%}

\end{keyword}

\end{frontmatter}

\newpage
\setcounter{equation}{0}
%\section{\label{sec:intro} Introduction}
It is well known that for a Poisson distribution with rate parameter $\lambda$, the mean and variance both equal $\lambda$
and the $n^{th}$ factorial moment is $\lambda^n$ for all $n\ge1$.
The Poisson distribution of order $k$ \cite{PhilippouGeorghiouPhilippou} is a variant or extension of the standard Poisson distribution.
For $k=1$ it is the standard Poisson distribution.
It is an example of a compound Poisson distribution \cite{Adelson1966}.
Its probability mass function is \cite{PhilippouGeorghiouPhilippou} 
\bq
P_n(k,\lambda) = e^{-k\lambda}\sum_{n_1+2n_2+\dots+kn_k=n} \prod_{j=1}^k \frac{\lambda^{n_j}}{n_j!} \,.
\eq
For $k=1$ this simplifies to standard result $P_n(\lambda) = e^{-\lambda}\lambda^n/n!$.
Philippou \cite{PhilippouMeanVar} derived the mean and variance
\begin{subequations}
\begin{align}
\mu_k &= (1+2+\dots+k)\lambda &=&\; \frac{k(k+1)}{2}\,\lambda \,,
\\
\sigma^2_k &= (1^2+2^2+\dots+k^2)\lambda &=&\; \frac{k(k+1)(2k+1)}{6}\,\lambda \,.
\end{align}
\end{subequations}
Less is known about its factorial moments.
Let $M_{(n)}(k,\lambda)$ denote the $n^{th}$ factorial moment of the Poisson distribution of order $k$.
The author obtained the following combinatorial sum 
\bq
\label{eq:facmom_orderk}
M_{(n)}(k,\lambda) = n!\sum_{n_1+2n_2+\dots+kn_k=n} \prod_{j=1}^k \frac{\lambda^{n_j}}{n_j!}\biggl[\binom{k+1}{j+1}\biggr]^{n_j} \,.
\eq
On searching the literature, it was discovered that 
Charalambides published the following expression for the factorial moment
(eq.~(4.9) in \cite{Charalambides1986})
\bq
\label{eq:Charalambides_eq4_9}  
M_{(n)} = n! \sum_{r=0}^n \lambda^r\,Q(n+r,r,k+1)/r! \,.
\eq
The quantity $Q(n+r,r,k+1)$ was given as follows (eq.~(3.19) in \cite{Charalambides1986})
\bq  
T_{n,r;k}(g_1,g_2,\dots,g_m) = n!k^{-r} Q(n+r,r,k+1)/r! \,.
\eq
Here (unnumbered before eq.~(3.19) in \cite{Charalambides1986})
\begin{subequations}
\begin{align}
\label{eq:gj}    
g_j &= \frac{j!}{k}\binom{k+1}{j+1} \,,
\\
\label{eq:Tnrk}    
\sum_{n=0}^\infty T_{n,r;k}(g_1,g_2,\dots,g_m)\frac{t^n}{n!} &= k^{-r}t^{-r} \biggl(\sum_{j=1}^k \binom{k+1}{j+1}t^{j+1}\biggr)^r \biggr/ r!\,.
\end{align}
\end{subequations}
({\em Note: this corrects a misprint in \cite{Charalambides1986}, where $t^{j+1}$ was written as $t^{i+1}$.})
The subscript ``$m$'' on $g_m$ is $m=\min\{n,k\}$.
Upon making the relevant substitutions into eq.~\eqref{eq:Charalambides_eq4_9},
the result coincides with eq.~\eqref{eq:facmom_orderk}.
Hence eq.~\eqref{eq:facmom_orderk} can be considered as a confirmation of the formulas in \cite{Charalambides1986}.
However \cite{Charalambides1986} did not display an explicit combinatorial sum 
for the factorial moments of the Poisson distribution of order $k$.

For the standard Poisson distribution ($k=1$),
the sum in eq.~\eqref{eq:facmom_orderk} contains only one summand and the result simplifies to (setting $n_1=n$)
\bq
\label{eq:facmom_orderk_k1}
M_{(n)}(1,\lambda) = n! \frac{\lambda^n}{n!}\biggl[\binom{2}{2}\biggr]^{n} = \lambda^n \,.
\eq
This is the well known result for the standard Poisson distribution.
However, for $k>1$, the sum in eq.~\eqref{eq:facmom_orderk} contains many more terms.
At a minimum, it is evident that $M_{(n)}(k,\lambda)$ is a polynomial in $\lambda$ of degree $n$.
It is also evident that $M_{(n)}(k,\lambda)=1$ for $n=0$ and has no constant term for $n>0$.
Also note that if $n<k$, then terms where $j>n$ do not appear in the sum in eq.~\eqref{eq:facmom_orderk} because we must have $n_j=0$ if $j>n$.
This tells us the following: for fixed $n\ge1$, the expression for $M_{(n)}(k,\lambda)$ is the {\em same} for all $k\ge n$.
For example:
\begin{enumerate}
\item
  For $n=1$, only $n_1$ is nonzero and the factorial moment is $M_{(1)}(k,\lambda) = \lambda$ for all $k\ge1$.
\item
  For $n=2$, only $n_1$ and $n_2$ are nonzero. We shall see below that 
\bq
\label{eq:M2}  
M_{(2)}(k,\lambda) = \biggl(\frac{k(k+1)}{2}\biggr)^2\lambda^2 + \frac{(k-1)k(k+1)}{3}\,\lambda \,.
\eq
This is the expression for $M_{(2)}(k,\lambda)$ for all $k\ge2$.
Observe that for $k=1$, the second term (in $\lambda$) vanishes and also $k(k+1)/2=1$ so the result is $M_{(2)}(1,\lambda)=\lambda^2$.
\end{enumerate}
For brevity of the notation in the formulas below, define the following (where $\kappa_j=0$ if $j>k$)
\bq
\kappa_j = \binom{k+1}{j+1} \,.
\eq
Hence $\kappa_1 = k(k+1)/2$ and $\kappa_2 = (k-1)k(k+1)/6$, etc.
Observe that $\kappa_j$  is proportional to $g_j$ in \cite{Charalambides1986} (see eq.~\eqref{eq:gj}).
Let us study the simplest case, viz.~$k=2$.
Set $n_2=s$ and $n_1 = n-2s$ and note that $0 \le s \le \lfloor(n/2)\rfloor$ to obtain
\bq
\label{eq:facmom_orderk_k2}
M_{(n)}(2,\lambda) = \sum_{s=0}^{\lfloor(n/2)\rfloor} \frac{n!}{(n-2s)!s!}\kappa_1^{n-2s}\kappa_2^s\,\lambda^{n-s} \,.
\eq
Then for $n=1,2,\dots$ (omitting explicit mention of the arguments $(k,\lambda)$)
\begin{subequations}
\begin{align}
M_{(1)} &= \kappa_1\lambda \,,
\\
M_{(2)} &= \kappa_1^2\lambda^2 +2\kappa_2\lambda \,,
\\
M_{(3)} &= \kappa_1^3\lambda^3 +6\kappa_1\kappa_2\lambda^2 \,,
\\
M_{(4)} &= \kappa_1^4\lambda^4 +12\kappa_1^2\kappa_2\lambda^3 +12\kappa_2^2\lambda^2 \,,
\\
M_{(5)} &= \kappa_1^5\lambda^5 +20\kappa_1^3\kappa_2\lambda^4 +60\kappa_1\kappa_2^2\lambda^3 \,,
\\
M_{(6)} &= \kappa_1^6\lambda^6 +30\kappa_1^4\kappa_2\lambda^5 +180\kappa_1^2\kappa_2^2\lambda^4 +120\kappa_2^3\lambda^3 \,.
\end{align}
\end{subequations}
Note the following:
\begin{enumerate}
\item
Clearly, for any $k>1$, even for $k=2$, the number of terms in $M_{(n)}$ is not bounded as a function of $n$.
For $k=2$, the polynomial $M_{(n)}(2,\lambda)$ contains all powers of $\lambda$ from $n$ down to $\lfloor(n+1)/2\rfloor$.
\item
More generally, the polynomial $M_{(n)}(k,\lambda)$ contains all powers of $\lambda$ from $n$ down to $\lfloor(n+k-1)/k\rfloor$.
It is only for $k=1$ that the number of terms in $M_{(n)}(k,\lambda)$ is a bounded function of $n$.
\item
  Stated another way, the lowest power of $\lambda$ in $M_{(n)}(k,\lambda)$ is $\lambda$ for $n\in[1,k]$
  and $\lambda^2$ for $n\in[k+1,2k]$ and in general it is $\lambda^j$ for $n\in[jk-k+1,jk]$.
\end{enumerate}

\newpage
Let us write the expressions for $M_{(n)}(k,\lambda)$ for low values of $n$ and arbitrary $k$.
We again omit explicit mention of the arguments $(k,\lambda)$.
\begin{subequations}
\begin{align}
M_{(1)} &= \kappa_1\lambda \,,
\\
\label{eq:M2_kappa}
M_{(2)} &= \kappa_1^2\lambda^2 +2\kappa_2\lambda \,,
\\
\label{eq:M3_kappa}
M_{(3)} &= \kappa_1^3\lambda^3 +6\kappa_1\kappa_2\lambda^2 +6\kappa_3\lambda \,,
\\
M_{(4)} &= \kappa_1^4\lambda^4 +12\kappa_1^2\kappa_2\lambda^3 +(12\kappa_2^2 +24\kappa_1\kappa_3)\lambda^2 +24\kappa_4\lambda \,,
\\
M_{(5)} &= \kappa_1^5\lambda^5 +20\kappa_1^3\kappa_2\lambda^4 +(60\kappa_1\kappa_2^2 +60\kappa_1^2\kappa_3)\lambda^3 +(120\kappa_1\kappa_4 +120\kappa_2\kappa_3)\lambda^2 +120\kappa_5\lambda \,,
\\
M_{(6)} &= \kappa_1^6\lambda^6 + 30\kappa_1^4\kappa_2\lambda^5
  +(180\kappa_1^2\kappa_2^2 +120\kappa_1^3\kappa_3)\lambda^4
  +(360\kappa_1^2\kappa_4 +120\kappa_2^3 +720\kappa_1\kappa_2\kappa_3)\lambda^3
  \nonumber\\
  &\quad\qquad\;\,
  +(720\kappa_1\kappa_5 +720\kappa_2\kappa_4 +360\kappa_3^2)\lambda^2 +720\kappa_6\lambda \,.
\end{align}
\end{subequations}
We adopt the notation in the textbook by Graham, Knuth and Patashnik \cite{GrahamKnuthPatashnik}
and write $[\lambda^j]M_{(n)}$ for the coefficient of $\lambda^j$ in $M_{(n)}$.
Then note the following:
\begin{enumerate}
\item
For all $k\ge1$ and $n\ge1$, the coefficient of the highest power $\lambda^n$ is
\bq
[\lambda^n]\,M_{(n)} = \frac{n!}{n!}\kappa_1^n = \biggl(\frac{k(k+1)}{2}\biggr)^n \,.
\eq
\item
  For all $k\ge2$ and $n\ge2$, the coefficient of the second highest power $\lambda^{n-1}$ is 
\bq
\begin{split}  
  [\lambda^{n-1}]\,M_{(n)} &= \frac{n!}{(n-2)!1!}\kappa_1^{n-2}\kappa_2
  \\
  &= \frac{n(n-1)}{3}\biggl(\frac{k(k+1)}{2}\biggr)^{n-1} (k-1) \,.
\end{split}  
\eq
\item
It is more complicated for the coefficient of the third highest power $\lambda^{n-2}$ and lower powers.
\item
  For all $k\ge3$ and $n\ge3$, the coefficient of the third highest power $\lambda^{n-2}$ is
\bq
\begin{split}
  [\lambda^{n-2}]\,M_{(n)} &= \frac{n!}{(n-4)!2!}\kappa_1^{n-4}\kappa_2^2 +\frac{n!}{(n-3)!1!}\kappa_1^{n-3}\kappa_3
  \\
  &= \frac16\binom{n}{3}\biggl(\frac{k(k+1)}{2}\biggr)^{n-2} (k-1)((2n-3)k -2n) \,.
\end{split}
\eq
\item
  For all $k\ge4$ and $n\ge4$, the coefficient of the fourth highest power $\lambda^{n-3}$ is 
\bq
\begin{split}
  [\lambda^{n-3}]\,M_{(n)} &= \frac{n!}{(n-6)!3!}\kappa_1^{n-6}\kappa_2^3
  +\frac{n!}{(n-5)!1!1!}\kappa_1^{n-5}\kappa_2\kappa_3
  +\frac{n!}{(n-4)!1!}\kappa_1^{n-4}\kappa_4
\\
&= \frac{2}{135}\binom{n}{4}\biggl(\frac{k(k+1)}{2}\biggr)^{n-3} (k-1)
   \biggl[(10n^2-45n+47)k^2 -5(4n^2-9n-1)k +2(5n^2+1)\biggr] \,.
\end{split}
\eq
\item
However, the coefficient of the {\em lowest} power $\lambda$ is simple, for all $k\ge1$ and $n\ge1$:
\bq
  [\lambda]\,M_{(n)} = \frac{n!}{1!}\kappa_n
  = n! \binom{k+1}{n+1} 
  = \frac{(k+1)k\dots(k-n+1)}{n+1} \,.
\eq
The value of the coefficient is zero if $k<n$.  
For $k\ge2$ and $n=2$, the coefficient of $\lambda$ is $(k-1)k(k+1)/3$, as stated in eq.~\eqref{eq:M2}.  
\end{enumerate}

\newpage
The factorial moments can also be derived from the factorial moment generating function (FMGF).
For the Poisson distribution of order $k$, the FMGF is, omitting the arguments $(k,\lambda)$
(\cite{Adelson1966}, see also \cite{Charalambides1986})
\bq
\label{qe:Poisson_orderk_FMGF}
M(t) = e^{-k\lambda}e^{\lambda (t+t^2+\dots+t^k)} \,.
\eq
The $n^{th}$ factorial moment is given by 
\bq
M_{(n)} = \frac{d^nM(t)}{dt^n}\biggr|_{t=1} \,.
\eq
The first few factorial moments are as follows.
\begin{enumerate}
\item
  The first derivative is
\bq
\frac{dM(t)}{dt} = \lambda e^{-k\lambda}e^{\lambda (t+t^2+\dots+t^k)} (1+2t+\dots+kt^{k-1}) \,.
\eq
\item
  Evaluate at $t=1$ to obtain the first factorial moment
\bq
M_{(1)} = \lambda (1+2+\dots+k) = \binom{k+1}{2}\lambda = \kappa_1\lambda \,.
\eq
\item
  The second derivative is
\bq
\begin{split}
\frac{d^2M(t)}{dt^2} &= \lambda e^{-k\lambda}\, \frac{d\ }{dt}\Bigl[\, e^{\lambda (t+t^2+\dots+t^k)} (1+2t+\dots+kt^{k-1}) \,\Bigr]
\\
&= \lambda e^{-k\lambda}\,e^{\lambda (t+t^2+\dots+t^k)}\biggl[\, \lambda(1+2t+\dots+kt^{k-1})^2 + 2!\sum_{j=2}^k\binom{j}{2}t^{j-2} \,\biggr] \,.
\end{split}
\eq
\item
  Evaluate at $t=1$ to obtain the second factorial moment
\bq
\begin{split}
M_{(2)} &= \lambda^2 \binom{k+1}{2}^2 +2!\lambda \binom{k+1}{3} 
\\
&= \kappa_1^2\lambda^2 +2\kappa_2\lambda 
\\
&= \biggl(\frac{k(k+1)}{2}\biggr)^2\lambda^2 +\frac{(k-1)k(k+1)}{3}\lambda \,.
\end{split}
\eq
This agrees with eqs.~\eqref{eq:M2} and \eqref{eq:M2_kappa}.
\item
  The third derivative is
\bq
\begin{split}
  \frac{d^3M(t)}{dt^2} &= \lambda e^{-k\lambda}\,
  \frac{d\ }{dt}\biggl\{\, e^{\lambda (t+t^2+\dots+t^k)}\biggl[\, \lambda(1+2t+\dots+kt^{k-1})^2 + 2!\sum_{j=2}^k\binom{j}{2}t^{j-2} \,\biggr]\biggr\}
\\
&= \lambda e^{-k\lambda}\,e^{\lambda (t+t^2+\dots+t^k)}\biggl\{\, \lambda^2(1+2t+\dots+kt^{k-1})^3
\\
&\qquad\qquad\qquad\qquad\qquad\quad
+ 3\times 2!\lambda\biggl(\sum_{j=2}^k\binom{j}{2}t^{j-2}\biggr)(1+2t+\dots+kt^{k-1})
\\
&\qquad\qquad\qquad\qquad\qquad\quad
+3!\sum_{j=2}^k\binom{j}{3}t^{j-3} \,\biggr\} \,.
\end{split}
\eq
\item
  Evaluate at $t=1$ to obtain the third factorial moment
\bq
\begin{split}
  M_{(3)} &= \lambda^3\binom{k+1}{2}^3 +6\lambda^2\binom{k+1}{2}\binom{k+1}{3} +3!\lambda \binom{k+1}{4} 
\\
&= \kappa_1^3\lambda^3 +6\kappa_1\kappa_2\lambda^2 +6\kappa_3\lambda 
\\
&= \biggl(\frac{k(k+1)}{2}\biggr)^3\lambda^3
+2\biggl(\frac{k(k+1)}{2}\biggr)^2(k-1)\lambda^2
+\frac{(k-2)(k-1)k(k+1)}{4}\lambda \,.
\end{split}
\eq
This agrees with eq.~\eqref{eq:M3_kappa} and is the same expression for all $k\ge3$.
For $k=1$ the last two terms vanish and for $k=2$ the last term vanishes.
\end{enumerate}

%\newpage
In conclusion, a formula for the factorial moments of the Poisson distribution of order $k$ was published in 1986 in \cite{Charalambides1986},
although an explicit combinatorial sum was not displayed.
The present note displays an explicit combinatorial sum for the factorial moments of the Poisson distribution of order $k$,
and is in agreement with the formula in \cite{Charalambides1986}.
Unlike the standard Poisson distribution (the case $k=1$), for $k>1$ the structure of the factorial moments is much more complicated.
Some properties of the factorial moments of the Poisson distribution of order $k$ were elucidated in this note.

%\newpage
\section*{\label{sec:ack}Acknowledgements}
The author thanks Professor Charalambides for clarifying the details of some of the formulas in \cite{Charalambides1986}.

%\newpage

\end{document}